\documentclass[12pt]{article}
\usepackage{pgf,tikz}
\usepackage{calc,graphicx, complexity}
\usepackage{graphics}
\everymath{\displaystyle}
\usepackage{graphicx}
\usepackage{color}
\usepackage{amssymb}
\usepackage{amsmath}
\newtheorem{prethm}{{\bf Theorem}}

\newenvironment{thm}{\begin{prethm}{\hspace{-0.5
				em}{\bf .}}}{\end{prethm}}
\newtheorem{prelemma}{{\bf Lemma}}

\newtheorem{preex}{{\bf Example}}

\newtheorem{preprop}{{\bf Proposition}}

\newenvironment{prop}{\begin{preprop}{\hspace{-0.5em}{\bf .}}}{\end{preprop}}
\newtheorem{precor}{{\bf Corollary}}

\newenvironment{cor}{\begin{precor}{\hspace{-0.5
				em}{\bf .}}}{\end{precor}}
\newtheorem{preremark}{{\bf Remark}}

\newtheorem{preprob}{{\bf Problem}}

\newtheorem{predefin}{{\bf Definition}}

\newtheorem{preconj}{{\bf Conjecture}}

\newtheorem{preprobb}{{\bf Problem}}

\newtheorem{prelem}{{\bf Theorem}}

\newtheorem{precla}{{\bf Claim}}

\newenvironment{proof}{{\bf Proof.}\rm }{\hfill{$\Box$}}

\newtheorem{presolution}{{\bf Solution.}}

\def\newpic#1{}
\def\qed{\ifhmode\unskip\nobreak\fi\quad\ifmmode\Box\else$\Box$\fi}

\title{\vspace{-2.2cm}\Large\bf\noindent $\mathcal{O}(VE)$ time algorithms for the Grundy (First-Fit) chromatic number of block graphs and graphs with sufficiently large girth}
\author{\large\bf Manouchehr Zaker\footnote{mzaker@iasbs.ac.ir}
\vspace{5mm}\\
Department of Mathematics,\\
Institute for Advanced Studies in Basic Sciences,\\
Zanjan 45137-66731, Iran\\
}

\date{}

\begin{document}
\maketitle

\begin{abstract}
\noindent The Grundy (or First-Fit) chromatic number of a graph $G=(V,E)$, denoted by $\Gamma(G)$ (or $\chi_{_{\sf FF}}(G)$), is the maximum number of colors used by a First-Fit (greedy) coloring of $G$. To determine $\Gamma(G)$ is $\NP$-complete for various classes of graphs. Also there exists a constant $c>0$ such that the Grundy number is hard to approximate within the ratio $c$. We first obtain an $\mathcal{O}(VE)$ algorithm to determine the Grundy number of block graphs i.e. graphs in which every biconnected component is complete subgraph. We prove that the Grundy number of a general graph $G$ with cut-vertices is upper bounded by the Grundy number of a block graph corresponding to $G$. This provides a reasonable upper bound for the Grundy number of graphs with cut-vertices. Next, define $\Delta_2(G)={\max}_{u\in G}~ {\max}_{v\in N(u):d(v)\leq d(u)} d(v)$. We obtain an $\mathcal{O}(VE)$ algorithm to determine $\Gamma(G)$ for graphs $G$ whose girth $g$ is at least $2\Delta_2(G)+1$. This algorithm provides a polynomial time approximation algorithm within ratio $\min \{1, (g+1)/(2\Delta_2(G)+2)\}$ for $\Gamma(G)$ of general graphs $G$ with girth $g$.
\end{abstract}



\section{Introduction}

\noindent All graphs in this paper are undirected without any loops and multiple edges. In a graph $G=(V,E)$, $\Delta(G)$ denotes the maximum degree of $G$. For a vertex $v\in V$, $N(v)$ denotes the set of neighbors of $v$ in $G$. Define the degree of $v$ as $d_G(v)=|N(v)|$. Denote by $d_G(u,v)$ the distance between $u$ and $v$ in $G$. For a subset $S\subseteq V(G)$, $G[S]$ denotes the subgraph of $G$ induced by the elements of $S$. A complete graph on $n$ vertices is denoted by $K_n$. The girth $g(G)$ of a graph $G$ on $n$ vertices and $m$ edges is the length of its shortest cycle (if exists) and can be obtained with time complexity $\mathcal{O}(nm)$. The size of a largest complete subgraph in $G$ is denoted by $\omega(G)$. A vertex $v$ is cut-vertex of $G$ if removing $v$ increases the number of connected components in $G$. By a block in a connected graph $G$, we mean any maximal $2$-connected subgraph of $G$. In addition, a subgraph isomorphic to $K_2$ which is bridge (i.e. cut-edge) in $G$ is also considered as block of $G$. Let $G$ be a graph and $C$ and ${\mathcal{B}}$ be the set of cut-vertices and blocks in $G$, respectively. Following \cite{W}, we define the block-cutpoint graph of $G$ as a bipartite graph $H=(C,{\mathcal{B}})$ in which one partite set is $C$, i.e. the set of cut-vertices and the partite set ${\mathcal{B}}$ has a vertex $b_i$ corresponding to each block $B_i$ in ${\mathcal{B}}$. We include $vb_i$ as an edge if and only if the block $B_i$ contains the cut-vertex $v$. The degree of a cut-vertex $v$ in $H$ is the number of connected components in $G\setminus v$. A proper vertex coloring of $G$ is an assignment of colors $1, 2, \ldots$ to the vertices of $G$ such that every two adjacent vertices receive distinct colors. By a color class we mean a subset of vertices having a same color. The smallest number of colors used in a proper coloring of $G$ is called the chromatic number of $G$ and denoted by $\chi(G)$. A strong result of Zuckerman asserts that to approximate $\chi(G)$ within a factor of $|V(G)|^{1-\epsilon}$ is $\NP$-hard for any fixed $\epsilon >0$ \cite{Zu}. But to obtain optimal or nearly optimal proper vertex colorings is an essential need in many applications. Hence we need to use and study efficient coloring heuristics. The Grundy (or First-Fit) coloring is a well-known and fast proper vertex coloring procedure defined as follow.

\noindent A Grundy-coloring of a graph $G$ is a proper vertex coloring of $G$ consisting of color classes say $C_1, \ldots, C_k$ such that for each $i < j$ each vertex in $C_j$ has a neighbor in $C_i$. Grundy-coloring is an off-line version of online First Fit coloring. Let the vertices of a graph $G$ be presented according to an ordering $\sigma: v_1, \ldots, v_n$. The First-Fit coloring assigns color $1$ the $v_1$ and then for any $i\geq 2$, assigns the smallest available color to $v_i$. The maximum number of colors used by the First-Fit coloring over all orderings $\sigma$ is called the First-Fit chromatic number and denoted by $\chi_{_{\sf FF}}(G)$. The Grundy number of a graph $G$, denoted by $\Gamma(G)$ is the maximum number of colors used in any Grundy-coloring of $G$. It can be easily shown that $\Gamma(G)=\chi_{_{\sf FF}}(G)$ \cite{Z2}. The literature is full of papers concerning the Grundy number and First-Fit coloring of graphs e.g. \cite{BFKS,CH,CK,CS,FGSS,GL,HHB,KPT,TWHZ,Z1,Z2,ZS}. Another important application of Grundy-coloring is in the online colorings of graphs, where the vertices are introduced according to an arbitrary or random orders \cite{AS,BGLP,GL,I,KT}. Clearly, $\Gamma(G)\leq \Delta(G)+1$. An improvement of the latter bound was obtained in \cite{Z21}. For any graph $G$ and $u\in V(G)$, define
$$\Delta(u)=\max \{d(v):v\in N(u),~d(v)\leq d(u)\},~~~ \Delta_2(G)=\max_{u\in V(G)} \Delta(u).$$
\noindent It was proved in \cite{Z21} that $\Gamma(G)\leq \Delta_2(G)+1$. Note that $\Delta_2(G)\leq \Delta(G)$ and $\Delta(G)-\Delta_2(G)$ may be arbitrarily large. The $\NP$-completeness of determining the Grundy number was proved for the complement of bipartite graphs in \cite{Z1} and \cite{Z2} and for bipartite graphs in \cite{HS}. Kortsarz in \cite{Ko} proved that there is a constant $c>1$ so that approximating the Grundy number within $c$ is not possible, unless $\NP \subseteq \RP$, where $\RP$ stands for the class of problems solvable by a randomized polynomial time algorithm. Only few families of graphs are known for which the Grundy number has polynomial time solution. Hedetniemi et al. \cite{HHB} obtained a linear algorithm for trees. Telle and Proskurowski proved that $\Gamma(G)$ can be determined in polynomial time for graphs with bounded tree-width \cite{TP}. Computational complexity of Grundy number for various classes of graphs is studied in \cite{BFKS}.

\noindent Let $G$ be a graph and $u\in V(G)$. By $\Gamma_G(u)$ we mean the maximum integer $j$ such that there exists a Grundy-coloring of $G$ in which the color of $u$ is $j$. Denote by $A_G(u)$ the set consisting of colors $j$ such that there exists a Grundy-coloring of $G$ in which $u$ has color $j$. It was proved in \cite{CS} that the Grundy number is continuous. In other words, for each $j$ with $\chi(G)\leq j \leq \Gamma(G)$ there exists a Grundy-coloring of $G$ with exactly $j$ colors. The following proposition asserts that $A_G(u)$ is continuous.

\begin{prop}
Let $G$ be a graph and $u\in G$. Then $A_G(u)=\{1, 2, \ldots, \Gamma_G(u)\}$.\label{contin}
\end{prop}

\noindent \begin{proof}
Write for simplicity $t=\Gamma_G(u)$. Vertex $u$ receives color $t$ in some Grundy-coloring of $G$ but receives no color more than $t$ in any Grundy-coloring of $G$. We show that for each $j$, $1\leq j \leq t-1$, there exists a Grundy-coloring of $G$ in which $u$ receives color $j$. Let $C_1, \ldots, C_t, \ldots, C_k$ be color classes in a Grundy-coloring of $G$ in which the color of $u$ is $t$, i.e. $u\in C_t$. Let $\sigma$ be an ordering of $V(G)$ by which the later Grundy-coloring has been obtained. If we recolor the vertices of $G$ in which the vertices of $C_t$ are colored firstly, then $u$ receives color $1$. We extend this partial Grundy-coloring to the whole graph. Vertex $u$ in the whole Grundy-coloring has color $1$. Suppose now that $j$ is fixed and arbitrary with $2\leq j \leq t-1$. Define $H=G[C_1 \cup \ldots \cup C_{j-1} \cup C_t]$. Let $\sigma'$ be an ordering of $H$ obtained by restricting $\sigma$ to the vertices of $H$. It is clear that in the Grundy-coloring of $(H, \sigma')$, vertex $u$ receives color $j$. The Grundy-coloring of $H$ extends to a Grundy-coloring of $G$ in which the color of $u$ is $j$. Hence, $j\in A_G(u)$, as desired.
\end{proof}

\noindent We also use the concept of System of Distinct Representatives (SDR). Let ${\bf A}=(A_1, \ldots, A_m)$ be a collection of subsets of a set $Y$. A system of distinct representatives for ${\bf A}$ is a set of distinct elements $a_1, \ldots, a_m$ in $Y$ such that $a_i\in A_i$. SDR is defined similarly if ${\bf A}$ is a multiset. In this paper by a list $L$ we mean a set whose elements are $1, 2, \ldots, t$ for some integer $t\geq 1$. We represent such a set $L$ as $L=\{1, 2, \ldots, t\}$. Given a collection of not-necessarily distinct lists $\mathcal{L}=\{L_1, \ldots, L_k\}$, a SDR $D$ for $\mathcal{L}$ is list-SDR if $D$ has form $D=\{1, \ldots, d\}$, for some $d\geq 1$.

\noindent To present our algorithms we need some knowledge about a fast sorting algorithm. Assume that for some integer $k\geq 1$, an array of $k$ integers $n_1, \ldots, n_k$ is given such that for each $i$, $n_i\leq k+1$. Using the counting sort (see Page 194 in \cite{CLRS}) we can sort these integers non-decreasingly with time complexity $\mathcal{O}(k)$.  Assume that lists $L_1, \ldots, L_k$ are such that $|L_i|\leq k$ for each $i$. Using the counting sort and consuming $\mathcal{O}(k)$ time steps, we can arrange the elements of this list as $L'_1, \ldots, L'_k$ such that $|L'_1| \leq |L'_2| \leq \cdots \leq |L'_k|$. As a bypass result, the multiplicity of each cardinality $|L_i|$ is also obtained.

\noindent {\bf The outline of the paper is as follows.} In Section $2$ we present a deterministic algorithm such that for a block graph $G$ with $m$ edges and $c$ cut-vertices, determines $\Gamma(G)$ and $\Gamma(w)$ for all cut-vertices $w$ of $G$ with time complexity $c\times {\mathcal{O}}(m)$ (Theorem \ref{poly-block}). This algorithm gives an upper bound for $\Gamma(G)$ of general graphs $G$ with cut-vertices (Proposition \ref{grundy-block}). An upper bound is also obtained for $\Gamma(G)$ of blocks graphs $G$ in terms of $\omega(G)$. In Section $3$ we prove in Theorem \ref{polyGrun} that $\Gamma(G)$ can be determined by an $\mathcal{O}(nm)$ algorithm if the girth $g$ of $G$ is at least $2\Delta_2(G)+1$. Also $\Gamma(G)\geq \lfloor (g+1)/2 \rfloor$ can be decided with a same complexity. Proposition \ref{approxGrun} provides an $\mathcal{O}(nm)$ approximation algorithm within ratio $\min \{1,(g+1)/(2\Delta_2(G)+2)$ for $\Gamma(G)$.

\section{An $\mathcal{O}(VE)$ algorithm for $\Gamma(G)$ of block graphs}

\noindent In this section we obtain an $\mathcal{O}(VE)$ algorithm for the Grundy number of block graphs. A graph $G$ is block graph if every block in $G$ is a complete subgraph. Online colorings of block graphs is an unexplored research area but there are interesting results for a special family of block graphs i.e. forests \cite{AS,BGLP,GL}. Let $G$ be a block graph and $w$ a cut-vertex in $G$. Our algorithm is based on a partition of $V(G)$ into subsets $F_1, \ldots, F_k$ with the following properties. Set $F_1$ consists of the non-cut-vertices in $G$ and $F_2$ is the set of non-cut-vertices in $G\setminus F_1$ except $w$, if $w$ is a non-cut-vertex in $G\setminus F_1$. The other subsets are defined similarly. For each $i\geq 3$, $F_i$ consists of non-cut-vertices in $G\setminus (F_1\cup \ldots \cup F_{i-1})$ except $w$ if $w$ is a non-cut-vertex in $G\setminus (F_1\cup \ldots \cup F_{i-1})$. Finally $F_k=\{w\}$. To following proposition states how to obtain this partition algorithmically.

\begin{prop}
Let $G$ be a block graph and $w$ a cut-vertex in $G$. Then the partition $\mathcal{F}=\{F_1, \ldots, F_k\}$ corresponding to $(G,w)$ can be obtained by an $\mathcal{O}(|E(G)|)$ algorithm.\label{partition}
\end{prop}

\noindent \begin{proof}
We obtain a Breadth First Search Tree $T^w$ in $G$ rooted at $w$ by consuming $\mathcal{O}(|E(G)|)$ time steps. Since $G$ is a block graph then BFS tree $T^w$ has the following property. Let $u$ be an arbitrary cut-vertex in $G$. For each neighbor $v$ of $u$ the following holds. Either each neighbor $z$ of $v$ is a neighbor of $u$ too (i.e. $u,v,z$ belong to a same block) or otherwise $v$ is a cut-vertex and the block containing $u$ and $v$ is different from the block containing $v$ and $z$. This fact implies that the non-cut-vertices in $G$ are the vertices of degree one in $T^w$ and vice versa. Recall that $F_1$ is the set of non-cut-vertices in $G$. It follows that $F_1$ is the set of vertices of degree one in $T^w$. Obviously $G\setminus F_1$ is a block graph and with a similar argument, the set of non-cut-vertices in $G\setminus F_1$ (except $w$) is the set of vertices of degree one in $T^w\setminus F_1$ (except $w$). It follows that for each $i\geq 3$, $F_i$ equals to the set of vertices of degree one in $T^w\setminus (F_1\cup \ldots \cup F_{i-1})$ except $w$ if $w$ is a vertex of degree one in the later graph. We conclude that the problem of obtaining the partition $\mathcal{F}$ in $G$ is reduced to obtain subsets $F_1, \ldots, F_k$ with their properties in $T^w$. When we construct the search tree $T^w$ the vertices of $T^w$ is partitioned into say $D_0, D_1, \ldots, D_t$, for some integer $t$, such that $D_j=\{v:d_T(w,v)=j\}$, for each $j$. Now we scan the vertices from $D_k$ down to $D_1$ and finally $D_0=\{w\}$. For each vertex $v\in D_k$, define $f(v)=1$. The subsets $F_i$ will be described by $f$-value of vertices. Assume that for some $j\geq 0$ the $f$-value of vertices in $D_{j+1}$ is determined. For each vertex $u\in D_j$, define
$$f(u)=\max \{f(z):z\in D_{j+1}, uz\in E(T^w)\}+1.$$
\noindent If the set under the above definition is empty then the maximum is taken as zero and hence $f(u)=1$. Obviously each vertex $v$ with $f(v)=1$ has not any neighbor in lower level and hence has degree one in $T^w$. The converse of this fact also holds. By repeating this method for $T^w\setminus F_1$, we obtain by induction that $F_j=\{v\in T^w: f(v)=j\}$, for each $j$. We omit the full details. To estimate the running time of the procedure, we note that the $D_j$ sets are obtained during the BFS implementation with complexity $\mathcal{O}(|E(G)|)$. For each vertex $v$, $f(v)$ is determined by $\mathcal{O}(d_G(v))$ comparisons. Hence, the total time complexity to obtain the partition $\mathcal{F}$ is $\mathcal{O}(|E(G)|)$.
\end{proof}

\noindent For a general graph $G$ and a vertex $v\in G$, define $\Gamma_G(v)$ (or shortly $\Gamma(v)$) as the maximum color $t$ such that vertex $v$ is colored $t$ in a First-Fit coloring of $G$. Clearly, $\Gamma(v)\leq d_G(v)+1$ and $\Gamma(G)={\max}_{v\in G} \Gamma(v)$. Recall that using the counting sort (\cite{CLRS}, page 194), every sequence $n_1, \ldots, n_k$ satisfying $n_i\leq k+1$, for each $i$, can be sorted in non-decreasing form with time complexity $\mathcal{O}(k)$.

\begin{prop}
In any block graph $G$, there exists a vertex $w$ such that $\Gamma(w)=\Gamma(G)$ and every such vertex $w$ is cut-vertex if $\Gamma(G)>\omega(G)$.\label{block1}
\end{prop}

\noindent \begin{proof}
Obviously, there exists a vertex $w$ such that $\Gamma(w)=\Gamma(G)$. Assume that $\Gamma(G)>\omega(G)$ and $w$ is not cut-vertex. Let $B$ be a block of $G$ containing $w$. Since $w$ is not cut-vertex then $|B|=d_G(w)+1$. Then $\Gamma(G)=\Gamma(w)\leq d_G(u)+1 = |B| \leq \omega(G)$. This contradicts $\Gamma(G)>\omega(G)$.
\end{proof}

\noindent The following GRUNDY-BLOCK$(G,w)$ algorithm receives a block graph $G$ and a vertex $w$ and determines $\Gamma(w)$. The algorithm assigns to every vertex $v$ a list of colors of the form $L(v)=\{1, 2, \ldots, t\}$ (shortly, $L(v)=12\ldots t$), where $t$ depends on $v$. Initially, it assigns lists to non-cut-vertices. The list for $w$ is specified at the final step of the algorithm. When the algorithm performs to determine a list to an arbitrary cut-vertex $u$, there is already a set of lists in the neighborhood of $u$. Denote this set of lists by $\mathcal{L}(u)$. A list $L(u)$ for $u$ is obtained by the execution of ASSIGN-LIST$(G,u,\mathcal{L}(u))$ (or shortly ASSIGN-LIST$(G,u)$) (below) with the lists $\mathcal{L}(u)$ as the input.

\noindent {\bf Name:} ASSIGN-LIST$(G,u,\mathcal{L}(u))$\\
\noindent {\bf Input:} Graph $G$, vertex $u\in V(G)$ and a collection $\mathcal{L}(u)=\{L_1, \ldots, L_k\}$ of not necessarily distinct lists of elements appearing in the neighborhood of $u$ (the neighbors of $u$ corresponding to the lists $\{L_1, \ldots, L_k\}$ are also registered in the input).\\
\noindent {\bf Output:} A list $L(u)$ for $u$
\begin{enumerate}
  \item Increasingly sort the lists $L_1, \ldots, L_k$ in terms of their lengths and obtain the sorted sequence $L'_1, \ldots, L'_k$ such that $|L'_1| \leq \cdots \leq |L'_k|$.
  \item Pick the element $1$ from $L'_1$. Let $i_2> 1$ be the smallest index such that $2\in L'_{i_2}$ and pick element $2$ from $L'_{i_2}$. Let $i_3> i_2$ be the smallest index such that $3\in L'_{i_3}$ and pick element $3$ from $L'_{i_3}$. Continue and let $t$ be the last element which is picked by this method. \% {\it Note that $t$ is the largest integer such that $\{1, 2, \ldots, t\}$ is a list-SDR for $\{L_1, L_2, \ldots, L_k\}$.}
  \item Assign the list $\{1, 2, \cdots, t, t+1\}$ to $u$ as ASSIGN-LIST$(G,u)$ and return the lists $\{L'_1, L'_{i_2}, \ldots, L'_{i_t}\}$. \%
\end{enumerate}

\begin{figure}
\hspace*{2cm}
\begin{tikzpicture}			
\draw[black, thick] (1,0)-- (3,0);
\draw[black, thick] (3,0)-- (3,-2);
\draw[black, thick] (3,-2)-- (1,-2);
\draw[black, thick] (1,0)-- (3,-2);
\draw[black, thick] (1,0)-- (1,-2);
\draw[black, thick] (3,0)-- (1,-2);
\filldraw [black] (1,-2) circle(2pt)
node [anchor=east]{$1234$};
\filldraw [black] (3,-2) circle(2pt)
node [anchor=west]{$123$};
\filldraw [black] (1,0) circle(2pt)
node [anchor=east]{$1$};
\filldraw [black] (3,0) circle(2pt)
node [anchor=west]{$12$};
\draw[black, thick] (8,0)-- (10,0);
\draw[black, thick] (10,0)-- (10,-2);
\draw[black, thick] (10,-2)-- (8,-2);
\draw[black, thick] (8,0)-- (10,-2);
\draw[black, thick] (8,0)-- (8,-2);
\draw[black, thick] (10,0)-- (8,-2);
\draw[black, thick] (6,-2)-- (8,-2);
\draw[black, thick] (6,-2)-- (6,-4);
\draw[black, thick] (6,-4)-- (8,-2);
\filldraw [black] (7.5,-2) circle(0pt)
node [anchor=south]{$1234$};
\filldraw [black] (10,-2) circle(2pt)
node [anchor=west]{$123$};
\filldraw [black] (8,0) circle(2pt)
node [anchor=east]{$1$};
\filldraw [black] (10,0) circle(2pt)
node [anchor=west]{$12$};
\filldraw [black] (6,-2) circle(2pt)
node [anchor=east]{$12$};
\filldraw [black] (6,-4) circle(2pt)
node [anchor=east]{$1$};
\filldraw [black] (1,-2) circle(0pt)
node [anchor=north]{$u$};	
\filldraw [black] (8,-2) circle(2pt)
node [anchor=north]{$u$};	
\end{tikzpicture}
\caption{Assignment of lists to $(K_4,u)$ (left) and ASSIGN-LIST$(G,u)$ assigns the list $1234$ to $u$ using the lists $1, 1, 12, 12, 123$ in the neighborhood of $u$ (right)}\label{pic1}
\end{figure}
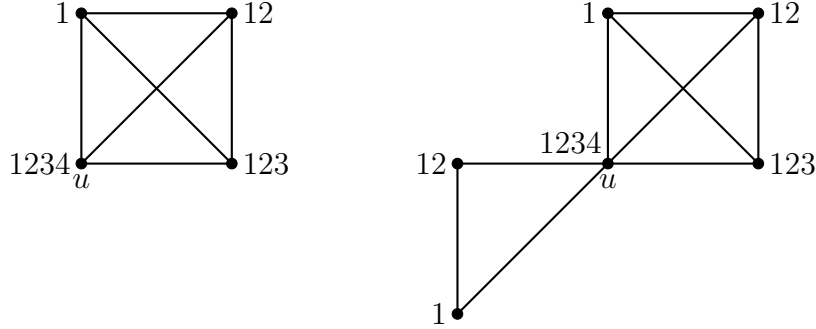

\noindent For example, for the lists $\{1\}, \{1\}, \{1,2\}, \{1,2\}, \{1,2,3,4,5\}$, ASSIGN-LIST$(G,u)$ outputs the list $\{1,2,3,4\}$. There are at most $k\leq d_G(u)$ list lengths in the command ${\bf 1}$, hence using the counting sort, the command ${\bf 1}$ takes ${\mathcal{O}}(d_G(u))$ time steps. The commands ${\bf 2}$ and ${\bf 3}$ take ${\mathcal{O}}(d_G(u))$ and ${\mathcal{O}}(1)$ time steps, respectively. Hence, for each $u$ the complexity of ASSIGN-LIST$(G,u,\mathcal{L}(u))$ is ${\mathcal{O}}(d_G(u))$. In case that $\mathcal{L}(u)=\emptyset$, ASSIGN-LIST$(G,u,\mathcal{L}(u))$ assigns $L(u)=\{1\}$. Therefore, if $B$ is a complete graph on say $q$ vertices then ASSIGN-LIST$(G,u)$ assigns the lists $1, 12, 123, \ldots, 12\cdots q$ to the vertices of $B$ as illustrated in Figure \ref{pic1} (left). ASSIGN-LIST$(G,u,\mathcal{L}(u))$ also returns $t$ neighbors of $u$ corresponding to the lists $\{L'_1, L'_{i_2}, \ldots, L'_{i_t}\}$. The following proposition is directly implied from these explanations.

\begin{prop}

\noindent $(i)$ For any collection $\mathcal{L}$ of not-necessarily distinct lists $L_1, \ldots, L_k$ such that for each $i$, $|L_i|\leq k$, ASSIGN-LIST outputs a list $\{1, \ldots, t\}$ within $\mathcal{O}(k)$ time steps so that $\{1, \ldots, t-1\}$ is a list-SDR for $\{L_1, \ldots, L_k\}$.

\noindent $(ii)$ ASSIGN-LIST$(G,u,\mathcal{L}(u))$ also obtains $t$ neighbors of $u$ say $u_1, \ldots, u_t$ such that $\{1, \ldots, t\}$ is a list-SDR for $\{L(u_1), \ldots, L(u_t)\}$ i.e. $i\in L(u_i)$, where $L(u)=\{1, \ldots, t\}$ is the output of the algorithm for the vertex $u$.\label{assign}
\end{prop}

\noindent In the following we introduce GRUNDY-BLOCK$(G,w)$. Except for the case $\Gamma(G)=\omega(G)$, in the light of Proposition \ref{block1}, to compute $\Gamma(G)$ we have to determine $\Gamma(w)$ for cut-vertices $w$. In fact, if $w$ is not cut-vertex then $\Gamma(w)=|B|$, where $B$ is a unique block containing $w$. In the following algorithm, every block graph is represented by its vertex partition $\mathcal{F}=\{F_1, \ldots, F_k\}$. Recall that $F_j$ is the set of non-cut-vertices in $G\setminus ({\bigcup}_{i=1}^{j-1} F_i)$ except $w$, if $w$ is a non-cut-vertex in the later set. Also $F_k=\{w\}$. By Proposition \ref{partition}, the partition $\mathcal{F}$ is obtained by an $\mathcal{O}(|E(G)|)$ algorithm.

\noindent {\bf Name:} GRUNDY-BLOCK$(G,w)$\\
\noindent {\bf Input:} A block graph $G$ with the associated partition $\{F_1, \ldots, F_k\}$ and a cut-vertex $w\in V(G)$\\
\noindent {\bf Output:} An assignment of lists to $V(G)$ in particular to $w$
\begin{enumerate}
  \item {\bf For} any vertex $v\in V(G)$, $\mathcal{L}(v) \leftarrow \emptyset$
  %\item $LV= \emptyset$ ~~~\% $LV$ is the set of vertices for which a list has be assigned.
  \item {\bf For} $i=1$ to $i=k$\\
  \hspace*{0.2cm} {\bf For} any $u\in F_i\setminus (F_1\cup \ldots \cup F_{i-1})$~~~ \% $F_0=\emptyset$ \\
\hspace*{0.4cm} Define $L(u)$ as ASSIGN-LIST$(G,u,\mathcal{L}(u))$\\
%\hspace*{0.6cm} $LV \leftarrow LV \setminus \{u\}$
 \hspace*{0.6cm} {\bf For} any neighbor $v\in N(u) \setminus (F_1\cup \ldots \cup F_{i-1})$,\\
\hspace*{0.9cm} $\mathcal{L}(v) \leftarrow \mathcal{L}(v) \cup \{L(u)\}$~~~ \% $\mathcal{L}(w)$ is updated at this step from its neighbors
  %\hspace*{1.4cm} $G \leftarrow G\setminus F_i$
  \item {\bf Return} $L(w)$ by ASSIGN-LIST$(G,w,\mathcal{L}(w))$.
\end{enumerate}

\noindent At the beginning, let $B$ be a block in $G$ with $t$ non-cut-vertices. Then GRUNDY-BLOCK$(G,w)$ assigns the lists $1, 12, \ldots, 12\cdots t$, to the non-cut-vertices of $B$. Note that if $v$ is an arbitrary non-cut-vertex in a block $B$ then GRUNDY-BLOCK$(G,w)$ can be performed such that the list of $v$ is $\{1\}$. This does not effect the output of the algorithm. These facts will be used in the proof of Proposition \ref{block}.

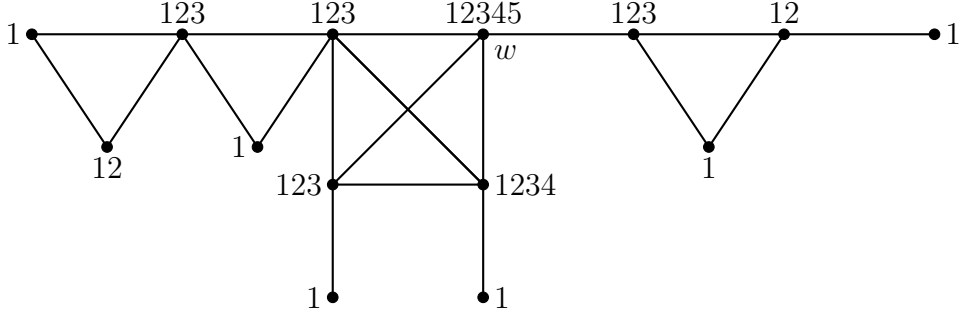
\begin{figure}
	\begin{center}
		\begin{tikzpicture}			

\draw[black, thick] (5,0)-- (3,0);
\draw[black, thick] (5,0)-- (7,-2);
\draw[black, thick] (7,0)-- (5,-2);
\draw[black, thick] (5,0)-- (7,-2);
\draw[black, thick] (5,0)-- (5,-2);
\draw[black, thick] (7,0)-- (7,-2);
\draw[black, thick] (5,-2)-- (7,-2);
\draw[black, thick] (9,0)-- (10,-1.5);
\draw[black, thick] (11,0)-- (10,-1.5);
\draw[black, thick] (5,-2)-- (5,-3.5);
\draw[black, thick] (7,-2)-- (7,-3.5);

\draw[black, thick] (4,-1.5)-- (5,0);
\draw[black, thick] (4,-1.5)-- (3,0);	
\draw[black, thick] (1,0)-- (3,0);
\draw[black, thick] (1,0)-- (2,-1.5);
\draw[black, thick] (5,0)-- (7,0);
\draw[black, thick] (7,0)-- (9,0);
\draw[black, thick] (9,0)-- (11,0);
\draw[black, thick] (11,0)-- (13,0);
\draw[black, thick] (2,-1.5)-- (3,0);	

\filldraw [black] (5,-2) circle(2pt)
node [anchor=east]{$123$};
\filldraw [black] (7,-2) circle(2pt)
node [anchor=west]{$1234$};
\filldraw [black] (10,-1.5) circle(2pt)
node [anchor=north]{$1$};
\filldraw [black] (5,0) circle(2pt)
node [anchor=south]{$123$};
\filldraw [black] (3,0) circle(2pt)
node [anchor=south]{$123$};
\filldraw [black] (2,-1.5) circle(2pt)
node [anchor=north]{$12$};
\filldraw [black] (1,0) circle(2pt)
node [anchor=east]{$1$};		
\filldraw [black] (7,0) circle(2pt)
node [anchor=south]{$12345$};		
\filldraw [black] (7.3,0) circle(0pt)
node [anchor=north]{$w$};	
\filldraw [black] (9,0) circle(2pt)
node [anchor=south]{$123$};		
\filldraw [black] (11,0) circle(2pt)
node [anchor=south]{$12$};		
\filldraw [black] (13,0) circle(2pt)
node [anchor=west]{$1$};		
\filldraw [black] (5,-3.5) circle(2pt)
node [anchor=east]{$1$};			
\filldraw [black] (7,-3.5) circle(2pt)
node [anchor=west]{$1$};
\filldraw [black] (4,-1.5) circle(2pt)
node [anchor=east]{$1$};
\end{tikzpicture}
\end{center}
\caption{A full execution of GRUNDY-BLOCK$(G,w)$}\label{pic2}
\end{figure}

\begin{prop}
Let $G$ be a block graph and $w$ a cut-vertex of $G$. For any $u\in V(G)$, let $L(u)$ be the list assigned by GRUNDY-BLOCK$(G,w)$ to $u$. Then $\Gamma(w)=|L(w)|$.\label{block}
\end{prop}

\noindent \begin{proof}
We first prove that $\Gamma(w)\leq |L(w)|$. Let $c$ be a Grundy-coloring of $G$ such that $c(w)=\Gamma(w)$. One of the following cases hold.

\noindent {\bf Case 1.} There exists a vertex $v$ of color $1$ in $c$ such that $v$ is a cut-vertex of $G$.

\noindent In this case, let $M$ be a connected component of $G\setminus \{v\}$ such that $w\not\in M$. Let $c'$ be the restriction of $c$ on $G'=G\setminus M$. Then we have $c'(w)=\Gamma(w)$. Denote by $L_{G'}(w)$ the output of GRUNDY-BLOCK$(G',w)$. Since $G'$ is induced subgraph of $G$ then $L_{G'}(w)\leq L(w)$. Also $|G'|<|G|$. By applying the induction hypothesis for $G'$ we obtain the following which proves the desired inequality in this case.
$$\Gamma(w)=\Gamma_{G'}(w)\leq L_{G'}(w) \leq L(w).$$

\noindent {\bf Case 2.} No vertex of color $1$ in $c$ is cut-vertex in $G$.

\noindent In this case, the algorithm GRUNDY-BLOCK$(G,w)$ can be performed such that every vertex of color $1$ receives list $\{1\}$ by the algorithm. The output $L(w)$ does not change. On the other hand, any vertex of color say $j\geq 2$ has a neighbor of color $1$ in $c$ and hence can not receive list $\{1\}$ by the algorithm. In fact their lists are at least $\{1,2\}$. It follows that
$$\big\{v\in V(G):c(v)=1\big\}=\big\{v\in V(G):L(v)=\{1\}\big\},~~C_1==\{v\in V(G):c(v)=1\}.$$
\noindent Let $H=G\setminus C_1$. Since $\Gamma(w)\geq 2$ then $w\not\in C_1$ and hence $w\in H$. For each vertex $u\in V(H)$, denote by $L'(u)$ the list obtained by GRUNDY-BLOCK$(H,w)$ to $u$. All vertices having list $\{1\}$ in $G$ are absent in $H$. Therefore every vertex of list $\{1,2\}$ in $G$ has now list $\{1\}$ in $H$. Let $v_i$ be the $i$-th vertex such that GRUNDY-BLOCK$(G,w)$ returns a list to $v_i$. Hence the algorithm assigns lists to the vertices of $G$ according to the order $v_1, v_2, \ldots, v_n$. Obviously $c(v_1)=1$ and $v_1\not\in H$. Hence, there exists $p$ with $1\leq p \leq n-1$ such that $V(H)=\{v_p, \ldots, v_n\}$. We prove by induction on $j\in \{p, \ldots, n\}$ that $|L'(v_j)|=|L(v_j)|-1$, for any $j$. Every vertex of list $\{1,2\}$ has now list $\{1\}$ in $H$. Then the assertion holds for the first element in $V(H)$. Assume that $v_i$ is an arbitrary vertex and for its all neighbors $u$ in $\{v_p, \ldots, v_{i-1}\}$ we have $|L'(u)|=|L(u)|-1$. GRUNDY-BLOCK$(G,w)$, obtains list $L(u)$ by only using the lists of previous neighbors of $u$. It follows that $|L'(v_i)|=|L(v_i)|-1$. We conclude that $|L'(w)|=|L(w)|-1$. We have now $\Gamma(w)=\Gamma_H(w)+1\leq |L'(w)|+1=|L(w)|$, as desired.

\noindent Now, we prove by induction on $|V(G)|$ that $\Gamma(w)\geq |L(w)|$. Consider the lists from GRUNDY-BLOCK$(G,w)$
and define $F=\{v\in V(G): L(v)=\{1\}\}$. Clearly, $F$ is a maximal independent set in $G$ and $w\not\in F$. Let $H=G\setminus F$ and $L_H(w)$ be the list obtained by GRUNDY-BLOCK$(H,w)$. By the induction $\Gamma_H(w)\geq |L_H(w)|$. We have also $|L(w)|\leq |L_H(w)|+1$ and $\Gamma_H(w)\leq \Gamma_G(w)-1$, since $F$ is maximal independent set in $G$. Note that $\Gamma_H(w)=\Gamma_G(w)-1$ does not necessarily hold. Combining the inequalities we obtain $\Gamma(w)\leq |L(w)|$. This completes the proof.
\end{proof}

\noindent The complexity of ASSIGN-LIST$(G,u)$ for every vertex $u$ is ${\mathcal{O}}(d_G(u))$. Corresponding to each vertex $u$, $d_G(u)$ actions are needed to update $\mathcal{L}(v)$ for neighbors $v$ of $u$. Also ${\mathcal{O}}(d_G(u))$ steps are executed to assign a list for $u$. It follows that the total time complexity of GRUNDY-BLOCK$(G,w)$ is
${\mathcal{O}}\big({\sum}_{u\in V(G)\setminus w} d_G(u)\big)={\mathcal{O}}(m)$.
\noindent We have $\Gamma(G)=\max \{\omega(G), \Gamma_1(G)\}$, where $\Gamma_1(G)=\max \Gamma_G(w)$, where the maximum is taken over all cut-vertices $w$ in $G$. The clique number $\omega(G)$ is the size of a maximum block in $G$ and then can be computed in ${\mathcal{O}}(|E(G)|)$ time steps. Also, $\Gamma_1(G)$ is determined by applying GRUNDY-BLOCK$(G,w)$ for all cut-vertices $w$ of $G$. Hence, the following is immediate.

\begin{thm}
There exists a deterministic algorithm such that for any block graph $G$ with $m$ edges and $c$ cut-vertices, determines $\Gamma(G)$ and $\Gamma(w)$ for all cut-vertices $w$ of $G$ with time complexity $c\times {\mathcal{O}}(m)$.\label{poly-block}
\end{thm}

\noindent In a graph $G$ denote by $\tilde{\Delta}(G)$ the maximum degree of a cut-vertex in the block-cutpoint graph $H=(C,{\mathcal{B}})$ of $G$. The degree of a cut-vertex $v$ in $H=(C,{\mathcal{B}})$ is the number of blocks of $G$ containing $v$. Let $G$ be an arbitrary graph and $u$ a vertex of $G$. In the following construction, by attaching a complete graph $K_p$ to $G$ at $u$ we mean a graph obtained by identifying vertex $u$ of $G$ with an arbitrary vertex in $K_p$.

\begin{prop}

\noindent (i) In a block graph $G$, let $\omega = \omega(G)$ and $\tilde{\Delta}=\tilde{\Delta}(G)$. Then $\Gamma(G)\leq \tilde{\Delta}(\omega-1)+1$.

\noindent (ii) For any $t\geq 2$ and $p\geq 2$, there exists a graph $G_{t,p}$ such that $\tilde{\Delta}(G)=t$, $\omega(G)=p$ and
$\Gamma(G)= t(p-1)+1=\tilde{\Delta}(\omega-1)+1$.\label{blockomega}
\end{prop}

\noindent \begin{proof}
\noindent To prove $(i)$, let $\Gamma(G)=k$ and $u$ be a cut-vertex such that $\Gamma(u)=k$. Let also $B_1, \ldots, B_t$ be all blocks in $G$ containing $u$. There are $k-1$ distinct colors appearing in the neighborhood of $u$ in $V(B_1)\cup \ldots \cup V(B_t)$. Suppose that $B_j$ contains $p_j$ distinct colors, for each $j=1, \ldots, t$. We have ${\sum}_{j=1}^t p_j = k-1$ and $|B_j|-1\geq p_j$. Hence, $\omega(G)-1\geq p_j$. It follows that $t(\omega-1) = {\sum}_{j=1}^t \geq k-1$. Finally
$$\Gamma(G)=k\leq t(\omega-1)+1 \leq \tilde{\Delta}(\omega-1)+1.$$
\noindent To prove $(ii)$, we use the fact that if $G$ is a general graph and $H$ is a graph obtained from $G$ by attaching a complete graph $K_p(u)$ (isomorphic to $K_p$) to each vertex $u$ of $G$ (these complete graphs are vertex disjoint) then $\Gamma(H)=\Gamma(G)+p-1$. To observe this fact, corresponding to each $u\in G$, we assign colors $1, \ldots, p-1$ to the vertices of $K_p(u)\setminus \{u\}$. Let $G_{1,p}$ be isomorphic to $K_p$. Let $G_{t,p}$ be obtained by attaching $K_p$ to each vertex of $G_{t-1,p}$. Obviously $\omega(G_{t,p})=p$. We may assume by induction on $t$ that $\Gamma(G_{t-1,p})=(t-1)(p-1)+1$. Now, $\Gamma(G_{t,p})= \Gamma(G_{t-1,p})+(p-1)=t(p-1)+1$, as desired.
\end{proof}

\noindent To obtain upper bounds for $\Gamma(G)$ in terms of $\omega(G)$ (such as Proposition \ref{blockomega}) is an important research tradition e.g. \cite{CK,GL,KPT,TWHZ,Z2}. Let $G$ be a general graph. Define a graph denoted by $G\uparrow {\mathcal{B}}$ as follows. Put an edge between every two non-adjacent vertices belonging to a same block of $G$. Obviously $G\uparrow {\mathcal{B}}$ is a block graph. By Proposition \ref{grundy-block} (below) $\Gamma(G)\leq \Gamma(G\uparrow {\mathcal{B}})$. The polynomial time algorithm in Theorem \ref{poly-block} determines $\Gamma(G\uparrow {\mathcal{B}})$. It follows that $\Gamma(G\uparrow {\mathcal{B}})$ is a polynomial time upper bound for $\Gamma(G)$.

\begin{prop}
For any graph $G$, $\Gamma(G)\leq \Gamma(G\uparrow {\mathcal{B}})$.
\label{grundy-block}
\end{prop}

\noindent \begin{proof}
We present a procedural proof. Let $c$ be a Grundy-coloring of $G$ using $k=\Gamma(G)$ colors. During the procedure we gradually produce a set $F\subseteq V(G)=V(G\uparrow {\mathcal{B}})$ and a coloring $c'$ using $k$ colors for the vertices of $F$ such that $c'=c$ on $F$ and $c'$ is a proper Grundy-coloring for the subgraph of $G\uparrow {\mathcal{B}}$ induced on $F$. This implies that $\Gamma(G\uparrow {\mathcal{B}}) \geq k$. At the beginning $F=\emptyset$. Let $u$ be a vertex of color $k$ in $c$. Pick a set $S_u$ consisting of $k-1$ neighbors of $u$ with different colors $1, 2, \ldots, k-1$. Add $u$ and $v$ to $F$ for any $v\in S_u$, define $c'(u)=c(u)$ and $c'(v)=c(v)$. Let $w$ be a neighbor of $u$ of color $k-1$ in $c$. We do the following for any color $i$ with $1\leq i \leq k-2$. If there exists a vertex $v\in F$ such that $c(v)=i$ and $w$ and $v$ belong to a same block then we do nothing concerning the pair $(w,i)$. In fact, in this case $w$ has a neighbor of color $i$ in $G\uparrow {\mathcal{B}}$ under the coloring $F$. Otherwise, let $v$ be a neighbor of $w$ of color $i$ in $c$. Define $c'(v)=c(v)$ and add $v$ in $F$.

\noindent Note that no two vertices belonging to a same block have a same color in $c'$. In other words, $c'$ is a proper coloring for $F$ (as an induced subgraph of $G\uparrow {\mathcal{B}}$) until so far. We repeat a same procedure for all vertices of color $k-2$ in $F$.
We scan the vertices of color $k-2$ in $F$ according to an arbitrary but fixed ordering. Let $t$ be any vertex of color $k-2$ in $F$. For each $i$ with $1\leq i \leq k-3$, if there exists a vertex $v\in F$ such that $c(v)=i$ and $t$ and $v$ belong to a same block then we do nothing concerning the pair $(t,i)$. Otherwise, let $x$ be a neighbor of $t$ of color $i$ in $c$. Define $c'(x)=c(x)$ and put $x$ in $F$. We repeat this procedure for other vertices of color $k-2$ in $F$. Update $F$ and note that $c'$ is proper for $F$ and for each $j\geq k-2$ and for each $i<j$ any vertex in $F$ of color $j$ has a neighbor in $F$ of color $i$ in the graph $G\uparrow {\mathcal{B}}$.

\noindent By repeating this technique for every other color $j$ with $j\leq k-3$ and all vertices in $F$ of color $j$, we finally obtain a subset $F$ of $V(G\uparrow {\mathcal{B}})$ and a Grundy-coloring $c'$ for $F$ using $k$ colors. The Grundy-coloring of $F$ is extended to a Grundy-coloring of whole $G\uparrow {\mathcal{B}}$ using at least $k$ colors. This completes the proof. The whole procedure uses ${\mathcal{O}}(|E(G)|)$ time steps.
\end{proof}

\noindent Proposition \ref{blockomega} yields a bound for the Grundy number of general graphs with cut-vertices. Let $\beta$ be the cardinality of a largest block in $G$. Obviously $\omega(G\uparrow {\mathcal{B}})=\beta$. Also the block-cutpoint graphs of $G$ and $G\uparrow {\mathcal{B}}$ are isomorphic and then the parameter $\tilde{\Delta}$ is same for both graphs. We obtain the final result of this section.

\begin{cor}
Let $G$ be a graph with at least one cut-vertex. Let $\beta$ be the size of maximum block in $G$ and $\tilde{\Delta}$ be the maximum degree of a cut-vertex of $G$ in its block-cutpoint graph. Then $\Gamma(G)\leq (\beta -1)\tilde{\Delta}+1$.
\end{cor}

\section{Graphs with sufficiently large girth}

\noindent In this section we obtain an exact polynomial time algorithm for $\Gamma(G)$ when girth of $G$ is sufficiently large with respect to a degree-related parameter (Theorem \ref{polyGrun}). This algorithm results in an approximation algorithm for the Grundy number of general graphs (Proposition \ref{approxGrun}). Given a graph $G$, a vertex $u\in V(G)$ and an integer $r\geq 0$, define a ball of radius $r$ centered at $u$ as $B(u,r)=\{v\in V(G):d_G(u,v)\leq r\}$, where $d_G(u,v)$ is the distance between $u$ and $v$ in $G$. We need to introduce a special subgraph in $G$. Let $u\in V(G)$, recall that $\Delta(u)=\max \{d(v):v\in N(u),d(v)\leq d(u)\}$. Define $G(u)=G[B(u,\Delta(u))]$. Note that $B\big(u,\Delta(u)\big)\subseteq B\big(u,d_G(u)\big)$. The following result proves that the Grundy number has a locality property.

\begin{prop}

\noindent (i) Let $G$ be a graph, $c$ be a Grundy-coloring of $G$ using $k$ colors. Let $u$ be a vertex in $G$ such that $c(u)=k$. Let also $H$ be a subgraph of $G$ with minimum cardinality such that $u\in H$ and the restriction of $c$ on $V(H)$ is a Grundy-coloring of $H$ with $k$ colors. Then $V(H)\subseteq B(u, k-1) \subseteq B(u, d_G(u))$. Also either $V(H)\subseteq B(u, \Delta(u))$ and then $V(H)\subseteq G(u)$ or for some neighbor $w$ of $u$, $V(H)\subseteq B(w, \Delta(w))$ and then $V(H)\subseteq G(w)$.

\noindent (ii) For any graph $G$, $\Gamma(G)=\max_{v\in V(G)} \Gamma(G(v))=\max_{v\in V(G)} \Gamma_{G(v)}(v)$.\label{Grunlocal}
\end{prop}

\noindent \begin{proof}
Let $C_1, \ldots, C_k$ be color classes in $c$. By the minimality of $H$, $V(H)\cap C_k =\{u\}$ and also $|V(H)\cap C_{k-1}|=1$. Let $w$ be an only neighbor of $u$ in $H$ such that $c(w)=k-1$. Since $H$ is minimal, each vertex in $V(H)\cap C_{k-2}$ has a neighbor in $V(H)\cap C_{k-1}$ and hence $V(H)\cap C_{k-2}\subseteq B(u,2)$. In general each vertex in $V(H)\cap C_j$ has a neighbor in $V(H)\cap C_{j+1}$ and hence by induction on $j$ we obtain $V(H)\cap C_j\subseteq B(u,k-j)$. It follows that $V(H)\subseteq B(u, k-1)$ and also $V(H)\subseteq B(w, k-1)$. Note that since $c(u)=k$ then $k-1\leq d(u)$ and $k-1\leq d(w)$. There are two possibilities.

\noindent {\bf (1)} If $d(w)\leq d(u)$ then $d(w)\leq \Delta(u)$. We have $V(H)\subseteq B(u, k-1) \subseteq B(u, d(w) \subseteq B(u, \Delta(u))$.

\noindent {\bf (2)} If $d(u)\leq d(w)$ then $d(u)\leq \Delta(w)$. We have $V(H)\subseteq B(w, k-1) \subseteq B(w, d(u) \subseteq B(w, \Delta(w))$. This completes proof of $(i)$.

\noindent To prove Part $(ii)$, let $c$ be a Grundy-coloring using $\Gamma(G)=k$ colors and $u$ be a vertex of color $k$ in $c$. Let $w$ be a neighbor of $u$ satisfying Part $(i)$. If Case {\bf (1)} holds then $G(u)$ has a partial Grundy-coloring with $k$ colors in which the color of $u$ is $k$ and then $\Gamma(G(u))\geq \Gamma_{G(u)}(u)\geq k$. If Case {\bf (2)} holds then $\Gamma(G(w))\geq \Gamma_{G(w)}(w) \geq k$. It follows that ${\max}_{v\in V(G)} \Gamma(G(v)) \geq {\max}_{v\in V(G)} \Gamma_{G(v)}(v) \geq k$. The inverse inequality is obtained by the clear fact that $\Gamma(G)\geq \Gamma(G(v))\geq \Gamma_{G(v)}(v)$ for any $v\in V(G)$.
\end{proof}

\noindent The following is easily obtained using Proposition \ref{Grunlocal} and induction on $k$.

\begin{cor}
Let $G$ be a graph and $u\in V(G)$. Let $c$ be a Grundy-coloring of $G$ such that $c(u)=k$. Let $H$ be a minimal subgraph of $G$ such that $u\in V(H)$ and the restriction of $c$ on $H$ is a Grundy-coloring $H$ using $k$ colors. Then $V(H)\subseteq B(u,k-1)$ and at most one vertex of $H$ has distance $k-1$ from $u$ in $G$.\label{atmost}
\end{cor}

\noindent Let $T$ be an arbitrary tree and $u$ a non-cut-vertex in $T$. Since every tree is a block graph then by considering $G=T$ and $w=u$, GRUNDY-BLOCK$(G,w)$ can also be performed for $(T,u)$. Denote the restricted algorithm by GRUNDY-TREE$(T,u)$. We don't need to repeat the commands of GRUNDY-TREE$(T,u)$ but summary the whole process to implement GRUNDY-TREE$(T,u)$ and obtain lists $L(v)$ for all vertices $v$ in $T$ (including $u$). By a breadth first search starting at $u$ in $T$ we obtain sets $D_j=\{v\in T:d_T(u,v)=j\}$, $j=0, 1, \ldots, t$. Then we obtain the partition sets $F_1, \ldots, F_k$ using the $f$-value of the vertices described in Proposition \ref{partition}. Finally we execute GRUNDY-BLOCK for $T$ and non-cut-vertex $u$ and obtain the necessary lists, in particular $L(u)$.

%\noindent {\bf Name:} GRUNDY-TREE$(T,u)$\\
%\noindent {\bf Input:} A tree $T$, a vertex $u\in V(T)$ and corresponding levels $F_1, \ldots, F_k$\\
%\noindent {\bf Output:} A list $L=\{1, 2, \ldots, t\}$ for $u$
%\begin{enumerate}
 % \item {\bf For} any vertex $v\in V(T)$, $\mathcal{L}(v) \leftarrow \emptyset$
 % \item {\bf For} $i=1$ to $i=k$\\
 % \hspace*{0.2cm} {\bf For} any $v\in F_i\setminus (F_1\cup \ldots \cup F_{i-1} \cup \{u\})$~~~ \% $F_0=\emptyset$ \\
%\hspace*{0.4cm} Define $L(v)$ as ASSIGN-LIST$(T,v,\mathcal{L}(v))$\\
 %\hspace*{0.6cm} {\bf For} any neighbor $w\in N(v) \setminus (F_1\cup \ldots \cup F_{i-1})$,\\
%\hspace*{0.9cm} $\mathcal{L}(w) \leftarrow \mathcal{L}(w) \cup \{L(v)\}$~~~ \% $\mathcal{L}(u)$ is updated at this step from its neighbors
 % \item {\bf Return} $L(u)$ by ASSIGN-LIST$(T,u,\mathcal{L}(u))$.
%\end{enumerate}

\noindent Let $w$ be an arbitrary neighbor of $u$. In a top-down drawing of $T$, where $u$ is the top-most vertex, let $T_w$ be a branch of $T$ rooted at $w$. GRUNDY-TREE$(T,u)$ is based on the partition $\mathcal{F}=\{F_1, \ldots, F_k\}$ for $(T,u)$. The corresponding partition for $(T_w,w)$ is obtained by the restriction of $\mathcal{F}$ on the vertices of $(T_w,w)$. It follows that the list assignment in $(T,u)$ restricted on $V(T_w)$ is identical to the list assignment in $(T_w,w)$ obtained by GRUNDY-TREE$(T_w,w)$. These fact are used in the proof of next result.

\begin{prop}
Let $\{1, \ldots, t\}$ be the list output by GRUNDY-TREE$(T,u)$. Then $\Gamma_T(u)=t$ and a Grundy-coloring of $T$ in which $u$ receives color $t$ can be obtained by a polynomial time algorithm.\label{treealgo}
\end{prop}

\noindent \begin{proof}
We prove by induction on the number of vertices that there exists a Grundy-coloring $c$ in $T$ such that $c(u)=t$. The list $L(u)$ obtained by GRUNDY-TREE$(T,u)$ is in fact a list-SDR of $t-1$ lists appearing in the neighborhood of $u$ such as $L(u_1), L(u_2), \ldots , L(u_{t-1})$. By Proposition \ref{assign} $(ii)$, ASSIGN-LIST$(T,u)$ returns the neighbors $u_1, \ldots, u_{t-1}$. Note that $|L(u_j)|=j$ for each $j=1, \ldots t-1$. Applying the induction for $(T_{u_j},u_j)$ we obtain that there exists a Grundy-coloring $c_j$ of $T_{u_j}$ such that $c_j(u_j)=j$. The colorings $c_1, \ldots, c_{t-1}$ do not overlap and hence are consistent. Combining these colorings we obtain a Grundy-coloring say $c$ of $T$ such that $c(u)=t$. This implies $\Gamma_T(u)\geq t$.

\noindent Denote by $L(v)$ the list assigned by GRUNDY-TREE$(T,u)$ to each vertex $v$. We shortly write $L(v)\geq q$ whenever $L(v)=\{1, \ldots, p\}$ and $p\geq q$, where $q$ is an integer. We prove by induction on $c(u)$ that in every rooted tree $(T,u)$, $L(u)\geq c(u)$, where $c$ is a Grundy-coloring for $T$. Let $c$ be a Grundy-coloring of $T$ with $c(u)=k$. Then for each $j=1, \ldots, k-1$, there exists a neighbor $u_j$ of $u$ such that $c(u_j)=j$. Let $c_j$ be the restriction of $c$ on $T_j$. It follows by the induction that $L(u_j)\geq j$, for each $j=1, \ldots, k-1$. We conclude that $L(u)\geq k$ and then $L(u)\geq \Gamma_T(u)$. \end{proof}

\noindent Let $u$ be a vertex in a graph $G$ and suppose that the girth $g$ of $G$ satisfies $g\geq 2\Delta(u)+1$, where $\Delta(u)=\max \{d(v):v\in N(u),d(v)\leq d(u)\}$. We grow a BFS tree rooted at $u$ and of depth $\Delta(u)$. Denote this tree by $T_{u,\Delta(u)}$. In fact $V(T_{u,\Delta(u)})=\{w\in V(G):~d_G(u,w)\leq \Delta(u)\}$. It follows that $V(G(u))=V(T_{u,\Delta(u)})$. The BFS tree $T_{u,\Delta(u)}$ is obtained in time complexity $\mathcal{O}(|E(G)|)$. By applying GRUNDY-TREE to the tree $T=T_{u,\Delta(u)}$ and vertex $u$ we obtain $\Gamma_T(u)$. Recall that, given a graph $G$ and $u\in V(G)$, $A_G(u)$ is the set consisting of colors say $j$ such that there exists a Grundy-coloring of $G$ in which $u$ receives color $j$.

\begin{prop}
Let $G$ be a graph of girth at least $2\Delta_2(G)+1$ and $u\in V(G)$. Let $T=T_{u,\Delta(u)}$ and $L(u)$ be the output of GRUNDY-TREE$(T,u)$. Then $A_{G(u)}(u)=L(u)$. For any $j\in L(u)$, a Grundy-coloring of $G$ in which $u$ receives color $j$ is obtained by GRUNDY-TREE$(T,u)$. In particular, a Grundy-coloring of $G$ in which $u$ receives color $\Gamma_{G(u)}(u)$ is obtained in a polynomial time.\label{propthm}
\end{prop}

\noindent \begin{proof}
Let $L(u)=\{1, 2, \ldots, t\}$. In order to prove $A_{G(u)}(u)=L(u)$ it suffices to prove $\Gamma_{G(u)}(u)=t$, since by Proposition \ref{contin}, $A_{G(u)}(u)=\{1, \ldots, \Gamma_{G(u)}(u)\}$. Note that $t\leq \Delta(u)$. We use the fact that $V(G(u))=V(T_{u,\Delta(u)})$. Proposition \ref{treealgo} provides a Grundy-coloring $c$ of $T$ such that $c(u)=t$. Let $H$ be a minimal subgraph in $T$ described in Corollary \ref{atmost} and $c'$ be the restriction of $c$ on $H$. Since the girth of $G$ is at least $2\Delta_2(G)+1$ then no two vertices in $V(H)\subseteq V(T)$ are adjacent in $G(u)$ unless their distance from $u$ in $T$ is exactly $\Delta(u)$. By Corollary \ref{atmost}, at most one vertex of $H$ has distance $t-1$ from $u$ in $T$. It follows that no two vertices in $H$ of identical color has distance $\Delta(u)$ from $u$. Hence, $c'$ is a partial Grundy-coloring in $G(u)$ using $t$ colors and can be extended to a Grundy-coloring of whole $G(u)$. Then $\Gamma_{G(u)}(u)\geq t$. Conversely, let $\Gamma_{G(u)}(u)=k$. Hence, there exists a Grundy-coloring $c''$ in $G(u)$ such that $c''(u)=k$. Note that $k-1\leq \Delta(u)$. Quiet similar to the previous arguments, if $H$ is a minimal subgraph with respect to $c''$ with $u\in H$ then $H\subseteq B(u,k-1)\subseteq B(u,\Delta(u))=V(T)$. It follows that $c''$ is a Grundy-coloring in $T$ with $c''(u)=k$. Then $t=\Gamma_T(u)\geq k$, as desired. \end{proof}

\noindent It is not hard to prove the following theorem.

\begin{thm}

\noindent (i) The Grundy number of every graph $G$ on $n$ vertices, with $m$ edges and of girth at least $2\Delta_2(G)+1$ can be determined by an $\mathcal{O}(nm)$ algorithm.

\noindent (ii) For any graph $G$ and any integer $k\leq (g+1)/2$, it can be decided in $\mathcal{O}(nm)$ time complexity whether $\Gamma(G)\geq k$ and in particular $\Gamma(G)\geq \lfloor (g+1)/2 \rfloor$.\label{polyGrun}
\end{thm}

\noindent \begin{proof}
\noindent To prove $(i)$, by the assumption on $G$, for any vertex $u$ of $G$, $g\geq 2\Delta(u)+1$ and then by Proposition \ref{propthm} we can apply GRUNDY-TREE for $T$ and $u$ to obtain $\Gamma_{G(u)}(u)$. On the other hand, $\Gamma(G)=\max \{\Gamma_{G(u)}(u):u\in V(G)\}$. Let $u_0$ be a vertex such that $\Gamma(G)=\Gamma_{G(u_0)}(u_0)$. It follows that $\Gamma(G)$ can be obtained in time complexity $\mathcal{O}(nm)$. Also GRUNDY-TREE$(T, u_0)$ provides a Grundy-coloring of $G$ using $\Gamma(G)$ colors, where $T$ is a BFS tree rooted at $u_0$ in $G(u_0)$.

\noindent To prove $(ii)$, let $1\leq k\leq (g+1)/2$ be an integer. Then $g\geq 2k-1$. By Proposition \ref{Grunlocal}, $\Gamma(G)\geq k$ if and only if there exists a $u$ such that $\Gamma_{G(u)}(u)\geq k$ and $k-1\leq \Delta(u)$. On the other hand, for each vertex $u\in V(G)$ we can decide whether $G'(u)$ has a Grundy-coloring in which $u$ receives color $k$, where $G'=G[B(u, k-1)]$. It follows that $\Gamma(G)\geq k$ is decided in time complexity $\mathcal{O}(nm)$. The extreme case is $k=\lfloor (g+1)/2 \rfloor$ and hence $\Gamma(G)\geq \lfloor (g+1)/2 \rfloor$ is decided by the procedure.
\end{proof}

\noindent Kortsarz proved in \cite{Ko} that there exists a constant $c>0$ such that $\Gamma(G)$ cannot be approximated within factor $c$ unless $\NP \subseteq \RP$. The following immediate corollary provides a $(g+1)/ (2\Delta_2(G)+2)$-approximation algorithm.

\begin{prop}
There exists an $\mathcal{O}(nm)$ time approximation algorithm within ratio $\min \{1,\frac{g+1}{2\Delta_2(G)+2}\}$ for the Grundy number of graphs $G$ of girth $g$.\label{approxGrun}
\end{prop}

\noindent \begin{proof}
We design an approximation algorithm $\mathcal{A}$ of time complexity $\mathcal{O}(nm)$ as follows. We first obtain the girth $g$ of $G$ with an algorithm of time complexity $\mathcal{O}(nm)$ based on BFS trees. Then using the algorithm of time complexity $\mathcal{O}(nm)$ presented in the proof of Theorem \ref{polyGrun}, we first decide whether or not $\Gamma(G)\geq \lfloor (g+1)/2 \rfloor$.

\noindent {\bf Case 1.} $\Gamma(G)\geq \lfloor (g+1)/2 \rfloor$.

\noindent In this case, $\mathcal{A}$ simply outputs the value $\lfloor (g+1)/2 \rfloor$. The latter value approximates $\Gamma(G)$ within $(g+1)/(2(\Delta_2(G)+1)$ since
$$\frac{g+1}{2(\Delta_2(G)+1)}\Gamma(G) \leq \frac{g+1}{2} \leq \Gamma(G).$$

\noindent {\bf Case 2.} $\Gamma(G)\leq \lfloor (g+1)/2 \rfloor$.

\noindent In this case, using the algorithm presented in the proof of Corollary \ref{polyGrun} $(ii)$ we obtain the exact value of $\Gamma(G)$ and a maximum Grundy-coloring.
          \end{proof}

\end{document}